\documentclass[10pt,journal,compsoc]{IEEEtran}

\usepackage{amsmath}
\usepackage{amssymb}
\usepackage[dvipdfmx]{graphicx}
\usepackage{color}
\usepackage{dsfont}

\ifCLASSOPTIONcompsoc

\usepackage{cite}
\usepackage{url}

\usepackage{theorem}
\theorembodyfont{\rmfamily}

\newtheorem{definition}{Definition}

\newtheorem{theorem}{Theorem}
\newtheorem{lemma}{Lemma}

\DeclareMathOperator{\sfdet}{\mathsf{det}}

\DeclareMathOperator{\sfspan}{\mathsf{span}}

\DeclareMathOperator*{\sfmax}{\mathsf{max}}

\newcommand{\mat}[1]{\left[\: \begin{matrix} #1 \end{matrix} \:\right]}
\newcommand{\pmat}[1]{\left(\; \begin{matrix} #1 \end{matrix} \:\right)}
\newcommand{\spliteq}[1]{\begin{split} #1 \end{split}}
\newcommand{\simode}[1]{\left\{\:  \begin{aligned} #1 \end{aligned} \right.}

\begin{document}

\title{Higher-Order Allan Variance for Atomic Clocks of Arbitrary Order:
Mathematical Foundation}

\author{Takayuki Ishizaki$^1$, Taichi Ichimura$^1$,
Takahiro Kawaguchi$^2$, Yuichiro Yano$^3$, Yuko Hanado$^3$
\thanks{$^1$ Tokyo Institute of Technology, 2-12-1, Ookayama, Meguro, Tokyo, 152-8552, Japan}
\thanks{$^2$ Gunma University, 1-5-1, Tenjincho, Kiryu, Gunma, 376-8515, Japan}
\thanks{$^3$ National Institute of Information and Communications Technology, 4-2-1, Nukui-Kitamachi, Koganei, Tokyo, 184-8795 , Japan}
\thanks{Manuscript received April 19, 2005; revised August 26, 2015.}
}

\IEEEtitleabstractindextext{%
\begin{abstract}
In this paper, we perform a time-domain analysis of the higher-order Allan variance for atomic clock models of arbitrary order.
Adopting a standard atomic clock model where the time series of the clock reading deviation is expressed as a Wiener or integrated Wiener process, we define the higher-order Allan variance as the mean squared higher-order difference of the clock reading deviation.
The main results of this paper are threefold.
First, we prove that the higher-order difference operation of the clock reading deviation, which can be interpreted as a linear aggregation with binomial coefficients, is not only sufficient but also necessary for a resulting aggregated time series to be an independent and identically distributed Gaussian process.
Second, we derive a complete analytical expression of the higher-order Allan variance, which consists of both time-dependent and time-independent terms.
Third, we prove that the higher-order Allan variance is time-independent if and only if the order of difference operation is greater than or equal to the order of the atomic clock model.
\end{abstract}

}

\maketitle

\section{Introduction}

Atomic clocks use the quantum transitions of atoms and ions as their frequency reference to output extremely precise and accurate waves at a constant frequency \cite{takamoto2005optical,huntemann2016single,schwindt2015miniature,kitching2018chip}.
This characteristic of atomic clocks supports modern society in a wide range of fields where time and frequency are important\cite{camparo2007rubidium}.
Examples include communications, positioning, broadcasting, undersea exploration, finance, space exploration, and geophysics\cite{rentsch2023100,reid2016leveraging,bondarescu2012geophysical}.  


Frequency stability quantifies the performance of a clock, which must be evaluated statistically due to the stochastic nature of atomic clocks.
In order to assess frequency stability from time series data, a clock statistical parameter is estimated, assuming that the corresponding stochastic process has ``ergodicity," a property such that a time average is equal to an ensemble average \cite{lindsey1976theory}. 
Ergodicity also requires that the stochastic process be at least weakly stationary. 
Without such stationarity, frequency domain analysis based on the Wiener-Khinchin theorem, and existing measures of frequency stability such as the Allan variance and its variants \cite{allan2016historical}, the total variance \cite{greenhall1999total}, and the parabolic variance \cite{vernotte2015parabolic}, are invalid. 
Therefore, a rational estimation of the clock statistical parameter requires that the time series data of the clock reading deviation be made into a stationary process.

It is well known that the clock reading deviation data of atomic clocks is not stationary. 
In particular, it was recently shown that small atomic clocks have not only linear frequency drift, but also non-negligible nonlinear frequency drift \cite{gardner2016second}.
These drift components confirm that it is not rational to estimate a clock statistical parameter directly from non-stationary time series data, because the estimated value does not converge to the actual parameter even in the limit of increasing the number of data. 
Therefore, a rational estimate includes generating stationary data by appropriately transforming the obtained non-stationary data while preserving equivalent information about the estimated clock statistical parameter.

The first contribution of this paper proves the ``uniqueness'' of the linear aggregation that makes the clock reading deviation data stationary.
One possible transformation that can make non-stationary data stationary is to take weighted sums of the data divided into fixed intervals, i.e., to apply linear aggregation to the time series data. 
In this paper, we mathematically prove that the aggregation weights, such that the resulting aggregated time series is an independent and identically distributed Gaussian process (implying ergodicity), are essentially unique.
Interestingly, such weights can be found as binomial coefficients. 
This result is stated in Theorem~1.

Theorem~1 shows that the unique aggregation operation consists of binomial coefficients and can be interpreted as a higher-order difference operation, which supports the rationality of existing frequency stability measures in terms of ``necessity." 
For example, the Allan variance applies a second-order difference operation to define a clock statistical parameter, and the Hadamard variance applies a third-order difference operation. 
From these facts, it can be deduced from mathematical necessity that the existing measures are rational in the sense of making the clock reading deviation data stationary.
However, one should keep in mind that it is necessary to apply a higher-order difference operation according to the order of the atomic clock, i.e., the number of its differential equations.

The second contribution of this paper is the derivation of a complete analytical expression of the higher-order Allan variance, which 
 is a clock statistical parameter defined as the mean squared higher-order difference of the clock reading deviation.
The derived expression is a generalization of the Allan and Hadamard variances in \cite{zucca2005clock}. 
The definition of the higher-order Allan variance is consistent with the Allan and Hadamard variances in that the dependence of clock noise parameters on the averaging time is identical, and the coefficient for the white-frequency noise is normalized.
This result is stated in Theorem~2.

The higher-order Allan variance is expected to be effective for analyzing the order of atomic clocks based on the separation of noise components in long-term data.  
In addition, the analytical expression derived in this paper can be used to verify the reliability of numerical analysis, since it is independent of numerical errors in difference operations.
We note that especially for higher-order difference operations, numerical errors due to digit drop generally become pronounced leading to unreliabile numerical calculations.

The higher-order Allan variance, deduced from the necessity for stationarization, is essentially equivalent to the higher-order structure function of the clock reading deviation. 
The structure function is a classical concept introduced by Kolmogorov in 1941 for the analysis of turbulence \cite{kolmogorov1941local,kolmogorov1941dissipation}, and is still used as a powerful tool in this field \cite{arenas2006on}. 
It has also been shown that the structure function corresponds to a generalization of the Allan variance  \cite{lindsey1976theory,greenhall1982fundamental}. 
Theorem~2 in this paper derives a complete analytical expression of the higher-order structure function of the clock reading deviation as a function of clock noise parameters from the time-domain state-space model. 
We note that the evaluation of oscillator stability by structure functions in previous studies assumes that the power spectral density of the noise components is given in the frequency domain.
This is clearly different from the modeling and analysis in the time domain of this paper.

The third and final contribution of this paper is to derive a necessary and sufficient condition for the higher-order Allan variance to be independent of time.
In particular, we prove that the higher-order Allan variance is time-independent if and only if the order of the difference operation is equal to or larger than the order of the atomic clock. 
This result is stated in Theorem 3. 
Although structure functions for time series data of clock reading deviations can be made stationary by a higher-order difference operation, to the best of the authors' knowledge, there is no result that mathematically proves that a higher-order difference operation is ``necessary'' for the variance to be time-independent.

The remainder of this paper is organized as follows. 
In Section~\ref{sec:moddes}, we introduce a standard atomic clock model of arbitrary order, where the time series of the clock reading deviation is expressed as a Wiener or integrated Wiener process.
An equivalent discrete-time model is also introduced.
In Section~\ref{sec:stationary}, we mathematically prove that the higher-order difference operation of the clock reading deviation is necessary and sufficient for a discrete-time aggregated process to be an independent and identically distributed Gaussian process.
In Section~\ref{sec:highallan}, we derive the complete higher-order Allan variance analytical expression, based on the fact proven in Section~\ref{sec:stationary}, for atomic clock models of arbitrary order.
Furthermore, we analyze the relationship between the difference operation and atomic clock model orders with application examples provided at the end of the section.
Section~\ref{sec:conc} concludes this paper with a summary of the results and brief description of future work.

\medskip
\noindent \textbf{Notation}~ 
We denote the set of real numbers by $\mathbb{R}$, 
the $i$th canonical unit vector, i.e., the $i$th column of the identity matrix, by $e_i$,
the transpose of a matrix $A$ by $A^{\sf T}$, 
the subspace spanned by vectors $v_1,\ldots,v_n$ by $\sfspan \{v_1,\ldots,v_n\}$,
the ensemble mean of a random variable $X$ by $\mathbb{E}[X]$,
the Dirac delta function by $\delta(\cdot)$,
and the binomial coefficient by
\[
\pmat{n\\k} := \frac{n(n-1)\cdots(n-k+1)}{k!}.
\]
The equality ``$:=$" is used in the sense that a new symbol on the left-hand side is defined by known symbols on the right-hand side.

\section{Description of Atomic Clock Models}\label{sec:moddes}

\subsection{Continuous-Time Models}

First, we introduce a standard atomic clock model.
From experimental observations, it is known \cite{galleani2008tutorial,galleani2010time} that the time series of clock reading deviation relative to an ideal time coordinate can be mathematically modeled as a continuous-time stochastic process
\begin{equation}\label{eq:intmodel}
\spliteq{
    y(t) &= \sum_{i=1}^n \frac{c_{i}}{(i-1)!}t^{i-1} \\
    & + \sum_{i=1}^n \int_0^t \int_0^{t_1} \cdots 
    \int_0^{t_{i-1}} v_i(t_i) dt_i \cdots dt_2 dt_1
}
\end{equation}
where $c_i$ is a scalar constant corresponding to an initial value, $n$ denotes the order of the atomic clock model, and  
\[
v_i=\{v_i(t):t\geq 0\}
\]
denotes a continuous-time white Gaussian process such that
\begin{equation}
\spliteq{
    & \mathbb{E}\bigl[ v_i(t) \bigr] =0 \\
    & \mathbb{E}\bigl[ v_i(t)v_i(t+\tau) \bigr] = q_i^2 \delta(\tau) \\
    & \mathbb{E}\bigl[ v_i(t)v_j(t+\tau) \bigr] = 0
    ,\quad \forall \tau \in \mathbb{R}, \quad i\neq j.
    }
\end{equation}
The integral equation model in \eqref{eq:intmodel} can be represented by a system of stochastic differential equations. 
In particular, we express its state equation  as
\begin{equation}\label{eq:diffeqs}
    \simode{
    \dot{x}_1(t) & = x_2(t) + v_1(t) \\
    \dot{x}_2(t) & = x_3(t) + v_2(t) \\
    &~~\vdots \\
    \dot{x}_{n-1}(t) & = x_{n}(t) + v_{n-1}(t) \\
    \dot{x}_n(t) & =  v_n(t).
    }
\end{equation}
If the initial values are given such that
\[
    x_{i}(0) = c_{i},\quad 
    i= 1,2,\ldots, n,
\]
then the first state variable $x_1$ satisfies
\[
    x_1(t) = y(t),\quad \forall t\geq 0.
\]
Furthermore, we denote the stacked vectors by
\[
    x(t):=\mat{x_1(t) \\ \vdots\\ x_n(t)}
    ,\quad
    v(t):=\mat{v_1(t) \\ \vdots\\ v_n(t)}.
\]
Then, the stochastic differential equation system model in \eqref{eq:diffeqs} can be written in the matrix form as
\begin{equation}\label{eq:ssmodel}
   \Sigma : \simode{
    \dot{x}(t) &= A x(t) + v(t) \\
    y(t) &= C x(t)}
\end{equation}
where $A\in \mathbb{R}^{n\times n}$ and $C\in \mathbb{R}^{1\times n}$ are defined as
\begin{equation}\label{eq:ssmat}
    A:=
    \mat{
    0 & 1 & 0 & \cdots & 0 \\
    0 & 0 & 1 & \cdots & 0 \\
    \vdots & & \ddots & \ddots & \vdots \\
    \vdots & & & 0 & 1 \\
    0 & 0 & \cdots & \cdots & 0
    }
    ,\quad
    C := \mat{1 & 0 & \cdots & 0}.
\end{equation}
The system noise $v$ is a white Gaussian process such that
\[
    \mathbb{E}\bigl[ v(t) \bigr] =0
    ,\quad
    \mathbb{E}\bigl[ v(t) v^{\sf T}(t+\tau) \bigr] = 
    Q \delta(\tau)
    ,\quad \forall t\geq 0.
\]
where the covariance matrix is diagonal and given as 
\[
Q:= \mathsf{diag}(q_1^2,\ldots,q_n^2).
\]
In the rest of this paper, we assume that all noise variances $q_1^2,\ldots,q_n^2$ are nonzero without loss of generality.

\subsection{Equivalent Discrete-Time Models}

For the discussion below, we introduce the following standard notion of a discrete-time stochastic process \cite{grimmett2020probability}.

\begin{definition}
A discrete-time stochastic process 
\[
f=\{f[k]: k=0,1,\ldots \}
\]
is said to be an \textit{independent and identically distributed Gaussian process} (IID Gaussian process) if 
\[
\mathbb{E}\bigl[ f[k] f^{\sf T}[k'] \bigr] =0, \quad \forall k\neq k',
\]
and every $f[k]$ follows an identical Gaussian distribution.
\end{definition}

An important property of the IID Gaussian process is that the statistics about ensemble can be characterized by the statistics about time, i.e., the ergodicity is satisfied.
In particular, it follows from the the law of large numbers that
\[
\frac{1}{K} \sum_{k=0}^{K-1} f[k]
\xrightarrow{{\rm a.s.}}
\mathbb{E}\bigl[ f[k] \bigr]
\]
as $K$ goes to infinity for every time instant $k$, where ``a.s." indicates ``almost surely."
This property is fundamental to estimate statistical parameters, such as mean and variance, from an observed time series data, as pointed out in \cite{lindsey1976theory}.

Consider a discrete-time sequence $\{t_0,t_1,\ldots\}$ in an ideal time coordinate.
The interval between two adjacent times is assumed to be constant, or equivalently 
\begin{equation}\label{eq:distime}
t_k := \tau k
,\quad k = 0,1,\ldots,
\end{equation}
where $\tau$ denotes a sampling period.
The relation between $x(t_k)$  and $x(t_{k+1})$ in \eqref{eq:ssmodel}
can be formally expressed as 
\[
    x(t_{k+1}) = 
    e^{A \tau} x(t_k) +
    \int_{t_k}^{t_{k+1}} e^{A(t_{k+1}-t)}v(t)dt.
\]
Based on this expression, we can derive a discrete-time model that is equivalent to the continuous-time model in \eqref{eq:ssmodel} over the discrete time sequence as
\begin{equation}\label{eq:ssmodeld}
    \simode{
    x[k+1] &= A_{\tau} x[k] + v[k] \\
    y_{\tau}[k] &= C x[k]}
\end{equation}
where the discrete-time transition matrix $A_{\tau}$ is defined as
\begin{equation}\label{eq:Atau}
A_{\tau} := e^{A \tau},
\end{equation}
the discrete-time state $x$ and clock reading deviation $y_{\tau}$  as
\[
    x[k]:=x(t_k)
    ,\quad
    y_{\tau}[k]:= y (t_k).
\]
Furthermore, the system noise 
\[
v=\{v[k]:k=0,1,\ldots\}
\] 
is an IID Gaussian process such that every $v[k]$ follows the Gaussian distribution such that
\[
    \mathbb{E}\bigl[ v[k] \bigr] =0
    ,\quad
    \mathbb{E}\bigl[ v[k] v^{\sf T}[k] \bigr] = Q_{\tau},
    \quad \forall k=0,1,\ldots
\]
where the covariance matrix is given as
\begin{align}\label{eq:defQk}
    Q_{\tau} := 
    \int_{0}^{\tau}
	e^{At}   Q e^{A^{\sf T}t}    dt.
\end{align}

It should be noted that, because of the structure of $A$ in \eqref{eq:ssmat}, the power of $A_{\tau}$ can be represented as
\begin{equation}\label{eq:matAr}
A_{\tau}^k = \sum_{m=0}^{n-1}
    \frac{(\tau k)^m}{m!} A^m,
\end{equation}
which is proven by the fact that
\[
    A^m = 0,\quad \forall m \geq n.
\]
Furthermore, the characteristic polynomial of $A_{\tau}$ is 
\[
\sfdet(\lambda I -A_{\tau}) = (\lambda -1)^n.
\]
This can be proven by the fact that $A_{\tau}$ is an upper triangular matrix whose diagonal elements are all one.
Thus, the lowest order annihilator polynomial 
\begin{equation}\label{eq:CH}
(A_{\tau} -I)^n =0
\end{equation}
is obtained by Cayley-Hamilton theorem \cite{bernstein2009matrix}.

\section{Stationarizing Operation}\label{sec:stationary}

\subsection{Formulation}

\begin{figure}[t]
\centering
\includegraphics[width = .90\linewidth]{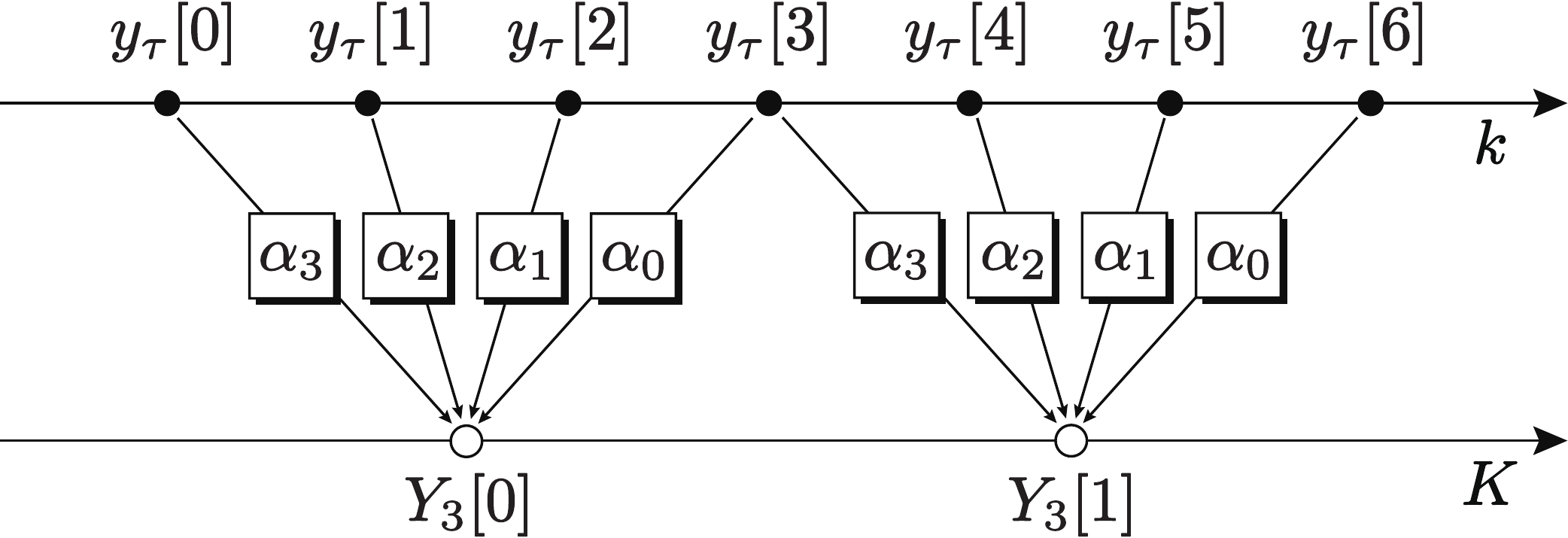}
\medskip
\caption{Construction of aggregated discrete-time process.}
\label{figYN}
\end{figure}

Consider the discrete-time stochastic process 
\[
y_{\tau}:=\{y_{\tau}[k]: k=1,2,\ldots\}.
\]
It is clear from the integral equation in \eqref{eq:intmodel} that this $y_{\tau}$ is not a stationary process because both mean and covariance are time dependent.
Therefore, in the following, we consider ``stationarizing" it by linear aggregation.
To this end, we introduce an aggregated discrete-time process as
\begin{equation} \label{eq:YNdef}
Y_N := \Biggl\{
\underbrace{
\sum_{i=0}^N \alpha_{i} y_{\tau}
\bigl[N(K+1)-i \bigr]
}_{Y_N[K]}
: K = 0,1,\ldots
\Biggr\}
\end{equation} 
where $\{\alpha_0,\ldots,\alpha_N\}$ is the set of nonzero weights for the linear aggregation, and $N$ specifies the number of signals to be aggregated.
This aggregation process is depicted in Fig.~\ref{figYN}, where $N$ is supposed to be 3.

By the definition of $Y_N[K]$ and the state equation of $\Sigma$ in \eqref{eq:ssmodel}, we see that
\begin{equation} \label{eq:YNKrep}
\spliteq{
Y_N[K] &=
\underbrace{
\left(\sum_{i=0}^{N}
\alpha_i
C A_{\tau}^{N-i} 
\right)
x[NK]
}_{X_N[K]}
\\
&+ 
\underbrace{
\sum_{i=0}^{N-1}
\sum_{j=0}^i
\alpha_j
C A_{\tau}^{i-j}
v[N(K+1)-1-i]
}_{V_N[K]}
.
}
\end{equation}
This expression will be used for the statistic analysis of the aggregated series $Y_N$.
It should be noted that the first and second terms are uncorrelated, namely
\begin{equation} \label{eq:XVun}
\mathbb{E}
\bigl[
X_N[K] V_N[K]
\bigr]
=0,
\quad
\forall K = 0,1,\ldots
\end{equation}
because those are, respectively, linear combinations of 
\[
\{v[0],\ldots, v[NK-1]\},
\quad
\{v[NK],\ldots, v[N(K+1)-1]\},
\]
which are independent each other.


\subsection{Analysis}

In this subsection, we analyze a condition that should be satisfied by the weight set for the aggregated process $Y_N$ in \eqref{eq:YNdef} to be an IID Gaussian process.
Such a weight set is essentially unique as shown in the following theorem.

\begin{theorem}\label{thm:IIDG}
For a an atomic clock model $\Sigma$ in \eqref{eq:ssmodel}, consider the discrete-time aggregated process $Y_N$ in \eqref{eq:YNdef}.
Then, there exists a weight set $\{\alpha_0,\ldots,\alpha_N\}$ such that $Y_N$ is an IID Gaussian process for any initial values of $\Sigma$ if and only if $N$ is greater than or equal to $n$.
In particular, if $N$ is equal to $n$, then such a weight set is uniquely determined except for constant multiples, and is found to be
\begin{equation} \label{eq:bicoe}
\alpha_i = \pmat{N\\i}(-1)^i,\quad i=0,\ldots,N.
\end{equation} 
\end{theorem}

\begin{IEEEproof}
We first prove that, if $N\geq n$, then $Y_N$ with the weight set in \eqref{eq:bicoe} is found to be an IID Gaussian process.
For $X_N[K]$ in \eqref{eq:YNKrep}, it follows that
\[
\sum_{i=0}^{N}
\pmat{
N \\ i
}
(-1)^{i}
A_{\tau}^{N-i}
= 
\underbrace{
(A_{\tau} - I)^n
}_{0} (A_{\tau} - I)^{N-n}=0,
\]
where \eqref{eq:CH} is used.
Furthermore,  $V_N[K]$ in \eqref{eq:YNKrep} is a linear combination of white Gaussian noises with zero means and constant variances.
Hence, $Y_N$ is proven to be an IID Gaussian process for any initial conditions.

Next, we prove that, if $Y_N$ for $N=n$ is supposed to be an IID Gaussian process, then the weight set is uniquely determined.
To this end, we first prove that, if $Y_n$ is an IID Gaussian process, then it follows that
\[
\beta := \sum_{i=0}^{n} \alpha_i CA_{\tau}^{n-i}
\]
is zero.
Denoting the $i$th element of $\beta$ by $\beta_i$, we have
\begin{equation} \label{eq:meanY}
\mathbb{E}\bigl[
Y_n[K]
\bigr]
=
\sum_{i=1}^n \beta_i \mathbb{E}\bigl[ x_i[nK] \bigr],
\end{equation} 
where the mean of $x_i[nK]$ is found to be
\[
\mathbb{E}\bigl[ x_i[nK] \bigr]
=
\sum_{m=0}^{n-i}  \frac{(\tau nK)^{m}}{m!} c_{i+m}.
\]
Note that only the mean of $x_n[nK]$ is constant.
This implies that, if any of $\beta_1,\ldots,\beta_{n-1}$ is not zero, then there exists some initial value such that \eqref{eq:meanY} depends on $K$.
Thus, all $\beta_1,\ldots,\beta_{n-1}$ must be zero if  $Y_n$ is an IID Gaussian process.
Furthermore, in the case of only $\beta_n$ is nonzero, we have
\[
\mathbb{E}
\left[
\left(
Y_n[K]
-\mathbb{E}
\bigl[
Y_n[K]
\bigr]
\right)^2
\right]
=
\beta_n^2 \left(
\mathbb{E}
\bigl[
x_n^2[nK]
\bigr]
-c_n^2
\right)
+
c
\]
where $c$ is a constant representing the variance of $V_n[K]$ in  \eqref{eq:YNKrep}.
This depends on $K$ because $x_n$ is a random walk process having a time-dependent variance.
This proves that $\beta$ must be zero if $Y_n$ is an IID Gaussian process.

In addition, we prove that, if $\beta$ is zero, then \eqref{eq:bicoe} is the unique weight set, except for constant multiples.
Using the expression in \eqref{eq:matAr}, we can rewrite $\beta$ as
\[
\beta =
\sum_{i=0}^n 
\sum_{m=0}^{n-1} 
\alpha_i \frac{ \{(n-i) \tau\}^m}{m!}
CA^m.
\]
Thus, the fact that $\beta$ is zero can be written as
\[
\underbrace{
\mat{1 & 1 & \cdots & 1 \\
0^1 & 1^1 & \cdots & n^1 \\
0^2 & 1^2 & \cdots & n^2 \\
\vdots & \vdots & \vdots & \vdots \\
0^{n-1} & 1^{n-1} & \cdots & n^{n-1}
}
}_{W_n}
\mat{
\alpha_n \\ \alpha_{n-1} \\ \alpha_{n-2} \\ \vdots \\ \alpha_0
}
=0.
\]
We  notice that $W_n \in \mathbb{R}^{n\times (n+1)}$ is a Vandermonde matrix \cite{bernstein2009matrix} being of row full rank.
Therefore, the kernel of $W_n$ is necessarily one-dimensional.
This proves that the weight set such that $\beta$ is zero is unique, except for constant multiples.

%
%

Finally, we prove that there exists a weight set such that $Y_N$ is an IID Gaussian process only if $N\geq n$.
This is equivalent to prove that, if $N<n$, then, for every nonzero weight set, $Y_N$ is not an IID Gaussian process for some initial value.
Analyzing a resultant Vandermonde matrix, we can verify that there does not exist a weight set such that
\[
\sum_{i=0}^{N} \alpha_i CA_{\tau}^{N-i} =0.
\]
Therefore, at least either of the mean or variance of $Y_N[K]$ depends on $K$.
This proves the claim.
\end{IEEEproof}
\medskip

Theorem~\ref{thm:IIDG} shows that there exists an essentially unique weight set for the aggregated process $Y_N$ in \eqref{eq:YNdef} to be an IID Gaussian process.
Furthermore, the number $N$ for the aggregation must be greater than or equal to the order $n$ of the atomic clock model.
The major novelty of this theorem is to show that, starting from the discussion of  ``necessity"  for stationarizing the time series of the sampled clock reading deviation, the linear aggregation to apply is uniquely found as the binomial coefficient in \eqref{eq:bicoe}.

\section{Higher-Order Allan Variance}\label{sec:highallan}

\subsection{Definition}

In this subsection, we introduce a generalized notion of the Allan variance for higher-order atomic clocks, based on the fact that the linear aggregation with the weight set in \eqref{eq:bicoe} can be interpreted in terms of a ``higher-order difference operation."
To this end, let us introduce the operation of difference for a continuous time series $f(\cdot)$.
In particular, the first order difference of $f(\cdot)$ is denoted as
\[
\triangle_{\tau} f(t) := f(t+\tau)- f(t).
\]
Then, the $N$th order difference is recursively defined as
\begin{equation}\label{eq:highdiff}
\triangle_{\tau}^{N} f(t) := 
\triangle_{\tau}  
\left\{
\triangle_{\tau}^{N-1}  
f(t)
\right\}.
\end{equation}
Using this difference operation, $Y_N[K]$ in \eqref{eq:YNdef} with the weight set in \eqref{eq:bicoe} can be represented as
\[
Y_N[K] = \triangle_{\tau}^{N} y(t_{NK}).
\]
From this viewpoint, Theorem~\ref{thm:IIDG} states that the higher-order difference of the sampled clock reading deviation is found to be an IID Gaussian process, and conversely, the higher-order difference must be used for the aggregation of the sampled clock reading deviation to be an IID Gaussian process.

Based on the notion of this higher-order difference, we introduce the following generalization of the Allan variance.

\begin{definition}\label{def:hoaver}
Consider the atomic clock model $\Sigma$ in \eqref{eq:ssmodel}.
The \textit{higher-order Allan variance} of $\Sigma$ is defined by
\begin{equation}\label{eq:hoallan}
\sigma^2_N(\tau;t) :=
\frac{1}{r_0^{(N)} \tau^2} \mathbb{E}
\left[
\left\{
\triangle^N_{\tau} y(t)
\right\}^2
\right],
\end{equation}
where $\triangle^N_{\tau}$ is the $N$th order difference operator in \eqref{eq:highdiff}, and 
\[
r_0^{(N)}:=
\sum_{i=0}^{N-1}
\sum_{j_1=0}^{i}
\sum_{j_2=0}^{i}
(-1)^{j_1 +j_2} \pmat{N \\j_1} \pmat{N \\j_2}.
\]
\end{definition}

We remark that $y(t)$ represents the clock reading deviation relative to an ideal time coordinate in Definition~\ref{def:hoaver}.
The constant $r_0^{(N)}$ is introduced for normalization.
The higher-order Allan variance coincides with the Allan and Hadamard variances for $N=2$ and $N=3$, respectively.
In particular, the Allan variance is obtained as
\[
\sigma_{2}^2(\tau; t)= \frac{1}{2\tau^2} 
\mathbb{E}
\left[
\left\{
y(t+2\tau) - 2 y(t+\tau) + y(t)
\right\}^2
\right],
\]
and the Hadamard variance is obtained as
\[
\spliteq{
\sigma_{3}^2(\tau; t)= \frac{1}{6\tau^2} 
\mathbb{E}
\Bigl[
\bigl\{
y(t+3\tau) &- 3 y(t+2\tau) \\
& + 3 y(t+\tau) - y(t)
\bigr\}^2
\Bigr].
}
\]
Interestingly, the higher-order Allan variance, mathematically deduced from the necessity for stationarization, also coincides with the higher-order structure function of the clock reading deviation, except for the difference of constant multiples \cite{lindsey1976theory}.
This coincidence affirms the rationality of existing oscillator stability measures.
We further remark that the normalization constant $r_0^{(N)}$ in \eqref{eq:hoallan}, which appears in the Allan and Hadamard variances, is first found by our time-domain analysis in Theorem~\ref{thm:highallan} below.

It has been pointed out in \cite{zucca2005clock} that the Allan variance of third order atomic clock models depends not only on the sampling period $\tau$, but also on the time $t$, while the Hadamard variance depends only on $\tau$.
This can be explained by the fact that the mean and variance of the second-order difference of the clock reading deviation $y(t)$ is time dependent while those of the third-order difference is not.
In the following, we aim at generalizing such a fact for atomic clock models of arbitrary order.

\subsection{Analytical Expression}

An analytical expression of the higher-order Allan variance can be derived as in the following theorem.

\begin{theorem}\label{thm:highallan}
Consider the atomic clock model $\Sigma$ in \eqref{eq:ssmodel}.
The higher-order Allan variance of $\Sigma$ is found to be
\begin{equation}\label{eq:highallan}
\spliteq{
\sigma^2_N(\tau;t)  = &
\sum_{i=1}^n \sum_{j=1}^n \frac{\gamma_i^{(N)} \gamma_j^{(N)}}{r_0^{(N)}}
\mathbb{E}\bigl[
x_i(t) x_j (t)
\bigr] 
\\
+ &
\sum_{m=0}^{n-1}
\frac{r_m^{(N)}}{r_0^{(N)}}
q_{1+m}^2
\tau^{2m-1}
}
\end{equation}
where $\gamma_m^{(N)}$ is a constant defined as
\[
\gamma_m^{(N)} := 
\frac{\tau^{m-2}}{(m-1)!} \sum_{i=0}^N 
\pmat{N\\i}(-1)^{N-i}
i^{m-1},
\]
and $r_m^{(N)}$ is a constant defined as
\[
\spliteq{
r_m^{(N)}  := & \left(\frac{1}{m!}\right)^2
\sum_{i=0}^{N-1}
\sum_{j_1=0}^{i}
\sum_{j_2=0}^{i}
(-1)^{j_1 +j_2} \pmat{N \\j_1} \pmat{N \\j_2}
 \\
\times & \sum_{l=0}^m
\pmat{m\\l}
\bigl\{
(i-j_1)(i-j_2)
\bigr\}^{m-l}
 \\
\times & \sum_{k=0}^l
\pmat{ l\\k }
\frac{(2i-j_1-j_2)^k}{2l-k+1}.
}
\]
\end{theorem}

\begin{IEEEproof}
We consider the case of $t=t_{NK}$.
This supposition does not lose generality because we can arbitrarily set the origin of the time variable such that $t$ is equal to $\tau NK$ with some $K$.
From \eqref{eq:YNKrep} and \eqref{eq:XVun}, we can see that
\[
r_0^{(N)}
\tau^2 \sigma^2_N(\tau;t_{NK})
=
\mathbb{E}\bigl[
X_N^2[K]
\bigr]
+
\mathbb{E}\bigl[
V_N^2[K]
\bigr].
\]
In the following, we evaluate these terms.
Using the expression in \eqref{eq:matAr} and the weight set in \eqref{eq:bicoe}, we have
\[
\spliteq{
X_N[K] & = 
\sum_{i=0}^{N} \alpha_i CA_{\tau}^{N-i} x(t_{NK}) \\
& = 
\sum_{i=0}^{N} \alpha_{N-i} CA_{\tau}^{i} x(t_{NK}) \\
&=\tau \sum_{m=1}^{n}  \gamma_m^{(N)} x_{m} (t_{NK}).
}
\]
Therefore, we obtain
\[
\mathbb{E}\bigl[
X_N^2[K]
\bigr]
= \tau^2 \sum_{i=1}^n \sum_{j=1}^n \gamma_i^{(N)} \gamma_j^{(N)} 
\mathbb{E}\bigl[
x_i (t_{NK}) x_j (t_{NK})
\bigr],
\]
which corresponds to the first term of \eqref{eq:highallan}.

We evaluete the mean of $V_N^2[K]$.
From \eqref{eq:YNKrep}, we have
\[
\spliteq{
\mathbb{E}\bigl[
V_N^2[K]
\bigr]
& =
 \sum_{i=0}^{N-1}
 \underbrace{
\Biggl(
\sum_{j=0}^i
\alpha_j
 C A_{\tau}^{i-j}
\Biggr)
}_{X_i}
Q_{\tau}
\Biggl(
\sum_{j=0}^i
\alpha_j
 C A_{\tau}^{i-j}
\Biggr)^{\sf T}
\\
&= 
\sum_{i=0}^{N-1}
\int^{\tau}_0
X_i e^{At} Q 
\left(X_i e^{At}
\right)^{\sf T} 
dt,
}
\]
where we have used the definition of  $Q_{\tau}$  in \eqref{eq:defQk}.
Note that
\[
\spliteq{
X_i e^{At} 
&=
\sum_{j=0}^i
\alpha_j
e_1^{\sf T}
e^{A\left\{(i-j)\tau +t \right\}} \\
& =
\sum_{j=0}^i
\alpha_j
\sum_{m=0}^{n-1}
\frac{\left\{ (i-j)\tau +t \right\}^m}{m!} e_{1+m}^{\sf T}.
}
\]
Thus, using the fact that
\[
e_i^{\sf T} Q e_j =
\left\{
\begin{array}{cl}
q_i^2, & i=j \\
0, & {\rm otherwise},
\end{array}
\right.
\]
we obtain
\[
\spliteq{
\mathbb{E}\bigl[
V_N^2[K]
\bigr]
& =
\sum_{i=0}^{N-1}
\sum_{m=0}^{n-1}
\int^{\tau}_0
\Biggl(
\sum_{j=0}^{i}
\alpha_j
\tfrac{\left\{ (i-j)\tau +t \right\}^m}{m!} q_{1+m}
\Biggr)^2
dt
\\
&=
\sum_{i=0}^{N-1}
\sum_{m=0}^{n-1}
\left(
\frac{q_{1+m}}{m!}
\right)^2
I_{im}
}
\]
where $I_{im}$ denotes the integral term of
\[
I_{im} :=
\int^{\tau}_0
\Biggl(
\sum_{j=0}^{i}
\alpha_j
\left\{ (i-j)\tau +t \right\}^m 
\Biggr)^2
dt.
\]
To evaluate the integral $I_{im}$, we use the expansions of
\[
\spliteq{
\Biggl(\sum_{j=0}^i x_j \Biggr)^2 & = \sum_{j_1=0}^i \sum_{j_2=0}^i x_{j_1} x_{j_2} , \\
(\alpha + \beta + \gamma)^m
& =
\sum_{l=0}^m
\pmat{m\\l}
\alpha^{m-l}
\sum_{k=0}^l
\pmat{l\\k}
\beta^{l-k}
\gamma^k.
}
\]
Applying these equalities, we can rewrite the integrand in $I_{im}$ just as a polynomial of $t$. 
Thus, calculating the integral, we have the second term of \eqref{eq:highallan}.
\end{IEEEproof}
\medskip

The analytical expression in Theorem~\ref{thm:highallan} is composed of two terms.
The first term is a time-dependent term that is related to the variance and covariance of atomic clock states.
The second term is a time-independent term that is represented as a polynomial of the sampling period $\tau$.

It is clear that whether the higher-order Allan variance is time dependent or not depends on whether the first term is nonzero or not.
In fact, an explicit expression of the time-dependent term can also be derived as follows.

\begin{table*}[ht]
\centering
\begin{tabular}{ccccccccccc}
$N$ & $q_1^2 \tau^{-1}$ & $q_2^2 \tau$ &  $q_3^2 \tau^{3}$ & $q_4^2 \tau^{5}$ 
 & $q_5^2 \tau^{7}$ & $q_6^2 \tau^{9}$ & $q_7^2 \tau^{11}$ & $q_8^2 \tau^{13}$  & $q_9^2 \tau^{15}$  & $q_{10}^2 \tau^{17}$ 
 \\ \hline
2 & 1.0000e+00 &  3.3333e-01 \\
3 & 1.0000e+00  & 1.6667e-01 &  9.1667e-02 \\ 
4 & 1.0000e+00 &  1.3333e-01 &  3.3333e-02 &  2.3968e-02 \\ 
5 & 1.0000e+00 &  1.1905e-01 &  2.2619e-02 &  6.9444e-03 &  6.1488e-03 \\ 
6 & 1.0000e+00 &  1.1111e-01 &  1.8254e-02 &  4.1005e-03 &  1.4863e-03 &  1.5632e-03 \\ 
7 & 1.0000e+00 &  1.0606e-01 &  1.5909e-02 &  3.0123e-03 &  7.7687e-04 &  3.2460e-04 &  3.9542e-04 \\ 
8 & 1.0000e+00 &  1.0256e-01 &  1.4452e-02 &  2.4531e-03 &  5.2278e-04 &  1.5218e-04 &  7.2018e-05 &  9.9720e-05  \\ 
9 &1.0000e+00 &  1.0000e-01 &  1.3462e-02 &  2.1170e-03 &  3.9850e-04 &  9.4365e-05 &  3.0604e-05 &  1.6180e-05 &  2.5098e-05 \\
10 & 1.0000e+00 &  9.8039e-02 &  1.2745e-02 &  1.8943e-03 &  3.2660e-04 &  6.7492e-05 &  1.7582e-05 &  6.2864e-06 &  3.6723e-06 &  6.3080e-06 \\
\hline
\end{tabular}
\medskip
\caption{Coefficients of higher-order Allan variance.}
\label{tab:coefnum}
\end{table*}

\begin{lemma}\label{lem:xcov}
Consider the atomic clock model $\Sigma$ in \eqref{eq:ssmodel}.
Then
\begin{equation}\label{eq:covxij}
\spliteq{
\mathbb{E}\bigl[
x_i(t) x_j (t)
\bigr] 
& =
\Biggl(
\sum_{m=0}^{n-i}\frac{t^m}{m!}c_{i+m}
\Biggr)
\Biggl(
\sum_{m=0}^{n-j}\frac{t^m}{m!}c_{j+m}
\Biggr)
\\
&+ 
\sum_{m=\sfmax\{i,j\}}^n \rho_{ij}^{(m)} q_m^2 t^{2m-i-j+1},
}
\end{equation}
where $\rho^{(m)}_{ij}$ is a positive constant defined as
\[
\rho^{(m)}_{ij} :=
\frac{1}{(m-i)!(m-j)!(2m-i-j+1)}.
\]
\end{lemma}

\begin{IEEEproof}
The state $x_i(t)$ can be formally written as
\[
x_i(t) = 
\underbrace{
e_i^{\sf T} e^{At}x(0)
}_{\overline{x}_i (t)} + 
\underbrace{
\int_0^{t} e_i^{\sf T} e^{A(t-\tau)} v(\tau) d \tau
}_{ \xi_i (t)}.
\]
Because the mean of $\xi_i (t)$ is zero, we have
\[
\mathbb{E}\bigl[
x_i(t) x_j (t)
\bigr] 
=
\overline{x}_i (t) \overline{x}_j (t)
+
\mathbb{E}\bigl[
\xi_i(t) \xi_j (t)
\bigr] .
\]
The first term here corresponds to the first term of \eqref{eq:covxij}.
On the other hand, the second term can be written as
\[
\mathbb{E}\bigl[
\xi_i(t) \xi_j (t)
\bigr]
= 
\int_0^{t} 
\underbrace{
e_i^{\sf T} e^{A(t-\tau)} 
}_{X_i (\tau)}
Q 
\left(e_j^{\sf T} e^{A(t-\tau)} \right)^{\sf T}
d\tau.
\]
Note that $X_i(\tau)$ can be represented as
\[
X_i(\tau) = \sum_{m=i}^{n} 
\frac{(t-\tau)^{m-i}}{(m-i)!} e_m^{\sf T}.
\]
Therefore, we have
\[
\spliteq{
\int_0^{t}
X_i(\tau) Q X_j^{\sf T}(\tau)
d \tau
& =
\int_0^{t}
X_{\sfmax \{i,j\}}(\tau) Q X_{\sfmax \{i,j\}}^{\sf T} (\tau)
d \tau
\\
& = 
\int_0^{t}
\textstyle
\sum_{m=\sfmax\{i,j\}}^n\limits
\tfrac{(t-\tau)^{2m-i-j}  }{(m-i)!(m-j)!}
q_m^2
d \tau.
}
\]
This integral corresponds to the second term of \eqref{eq:covxij}.
\end{IEEEproof}
\medskip

From Lemma~\ref{lem:xcov}, we see that the variance and covariance of atomic clock states are nonzero polynomials of $t$.
In particular, the first term of \eqref{eq:covxij} corresponds to the deterministic trend of atomic clock states, which is dependent on a part of the initial values $c_1,\ldots,c_n$, while the second term corresponds to their stochastic behavior, which is dependent on a part of  the noise variances $q_1^2,\ldots, q_n^2$.

Furthermore, the following lemma proves an essential fact to determine the time dependence of \eqref{eq:highallan}.

\begin{lemma}\label{lem:bicsum}
The constant $\gamma_m^{(N)}$ in \eqref{eq:highallan} is zero if and only if $m$ is less than or equal to $N$.
\end{lemma}

\begin{IEEEproof}
Let us first prove that 
\begin{equation}\label{eq:gam0}
\gamma_m^{(N)}=0
,\quad
\forall m=1,\ldots, N.
\end{equation}
From the binomial theorem, we have
\begin{equation}\label{eq:bithm}
(z-1)^N = \sum_{i=0}^N
\pmat{N\\i}(-1)^{N-i}
 z^{i}.
\end{equation}
Considering $z=1$, we see that  $\gamma_1^{(N)}$ is zero.
Furthermore, the differentiation of \eqref{eq:bithm} with respect to $z$ leads to
\[
N(z-1)^{N-1} = \sum_{i=0}^N  
\pmat{N\\i}(-1)^{N-i}
 i z^{i-1}.
\]
Considering $z=1$ again, we see that $\gamma_2^{(N)}$ is zero.
Repeating such a procedure proves \eqref{eq:gam0}.

Next, we prove that $\gamma_{m}^{(N)}$ is not zero for $m> N$.
Using the notation of \eqref{eq:bicoe}, we define
\[
w_m := \mat{
0^{m-1} \\ 1^{m-1} \\ \vdots \\ N^{m-1}
}
,\quad
\alpha:= \mat{ \alpha_N \\ \alpha_{N-1} \\ \vdots \\ \alpha_0}.
\]
Note that $\gamma_m^{(N)}$ is a constant multiple of $\alpha^{\sf T}w_m$.
Then, we see that \eqref{eq:gam0} is equivalent to
\[
\sfspan\{\alpha\}^{\perp} = \sfspan\{ w_1,\ldots,w_{N} \},
\]
where $\perp$ denotes the orthogonal complement.
Furthermore, $w_i$ and $w_j$ are linearly independent if $i \neq j$.
Thus
\[
w_m \notin \sfspan\{\alpha\}^{\perp}
,\quad
\forall m=N+1,N+2,\ldots,
\]
or equivalently
\[
\alpha^{\sf T} w_m \neq 0
,\quad
\forall m=N+1,N+2,\ldots.
\] 
Hence, the claim is proven.
\end{IEEEproof}
\medskip

Lemma~\ref{lem:bicsum}  in conjunction with Lemma~\ref{lem:xcov}
implies that the time-dependent term in \eqref{eq:highallan} disappears if and only if $N$ is greater than or equal to $n$, proving the following theorem.

\begin{theorem}
Consider the atomic clock model $\Sigma$ in \eqref{eq:ssmodel}.
The higher-order Allan variance of $\Sigma$ is time independent if and only if $N$ is greater than or equal to $n$.
\end{theorem}

\begin{IEEEproof}
The sufficiency is directly proven by Theorem~\ref{thm:highallan} and Lemma~\ref{lem:bicsum}.
This is because, if $N \geq n$, then
\[
\gamma^{(N)}_m =0
,\quad \forall m=1,\ldots,n.
\]
To prove the necessity, we suppose that $N < n$.
From the analytical expression in Lemma~\ref{lem:xcov}, we see that 
\[
\sigma^2_N(\tau;t)
= 
\tfrac{
\left( \gamma_{N+1}^{(N)}  \right)^2 
}{r_0^{(N)}}
\!
\rho_{(N+1)(N+1)}^{(n)} q_n^2
t^{2n-2N-1}
\! +  o(t^{2n-2N-2}),
\]
which is time dependent because $\gamma_{N+1}^{(N)}$ is not zero as shown in Lemma~\ref{lem:bicsum}.
Hence, the claim follows.
\end{IEEEproof}
\medskip

\subsection{Example}

We demonstrate application examples of Theorem~\ref{thm:highallan}.
Consider the third-order atomic clock model.
For the second-order difference, we have
\begin{subequations}
\begin{equation}\label{eq:suballan}
\sigma^2_2(\tau;t)
=
\frac{1}{2} (c_3^2+q_3^2 t)\tau^2
+
q_1^2 \tau^{-1}+ \frac{1 }{3} q_2^2 \tau + \frac{ 23 }{60} q_3^2 \tau^3.
\end{equation}
Furthermore, for the third-order difference, we have
\begin{equation}\label{eq:subhad}
\sigma^2_3(\tau;t)
=
q_1^2 \tau^{-1}+ \frac{1 }{6} q_2^2 \tau + \frac{ 11 }{120} q_3^2 \tau^3.
\end{equation}
\end{subequations}
These two values coincide with the Allan variance and the Hadamard variance shown in \cite{zucca2005clock}.

\begin{figure}[t]
\centering
\includegraphics[width = .95\linewidth]{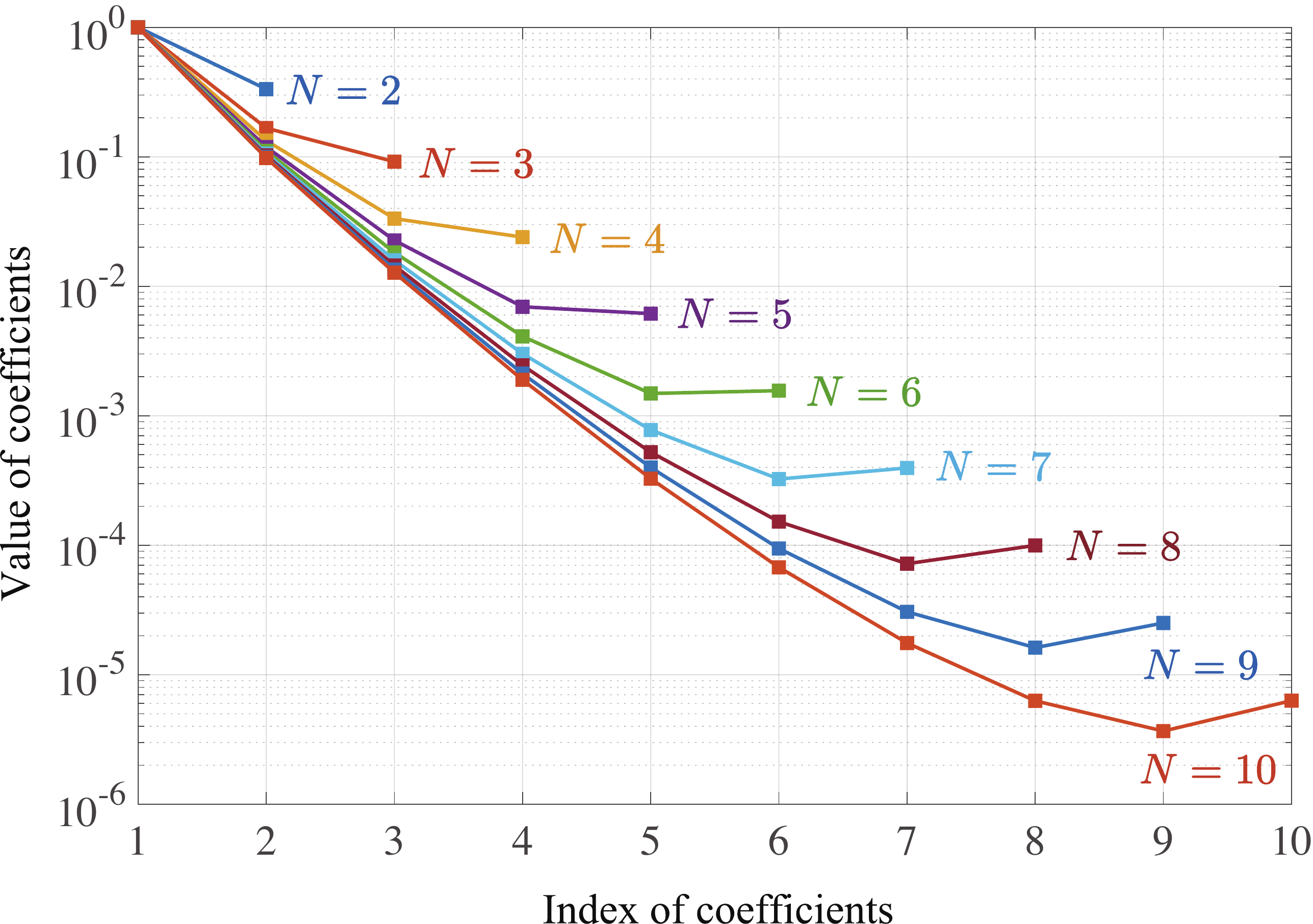}
\medskip
\caption{Logarithmic plot of coefficients of higher-order Allan variance.}
\label{figcoef}
\end{figure}

For reference, we numerically calculate the coefficients of the time-independent term of \eqref{eq:highallan} varying the difference order $N$.
We suppose that the order of atomic clock models is also $N$.
In this case, the higher-order Allan variance is not time dependent, and is simplified as
\[
\sigma_N^2(\tau) = \sum_{m=0}^{N-1} \frac{r_m^{(N)}}{r_0^{(N)}} q_{1+m}^2 \tau^{2m-1}.
\]
The result is shown in Tab.~\ref{tab:coefnum}, whose logarithmic plot is given in Fig.~\ref{figcoef}.
This result shows that coefficients with respect to noise variances for higher-order atomic clock models are exponentially smaller.
This trend is consistent with the empirical fact that longer term measurements are required to assess the accuracy of higher-order atomic clock models.

\section{Concluding Remarks}\label{sec:conc}

In this paper, we have mathematically deduced the definition of the higher-order Allan variance in terms of ``necessity" to make the stochastic process of clock reading deviation stationary.
In particular, we have proven that the aggregation weights are essentially unique for the resulting clock reading deviation to be an independent and identically distributed Gaussian process.
Furthermore, we have derived a complete analytical expression of the higher-order Allan variance for a standard atomic clock model of arbitrary order.

The main focus of this paper is to analyze a rational statistical parameter for the atomic clock model of arbitrary order.
Building an efficient algorithm to estimate the statistical parameter and extending the results to the overlapping version of the higher-order Allan variance is part of future work.



\section*{Acknowledgments}

This paper includes the results of research and development conducted by the Ministry of Internal Affairs and Communications (MIC) under its "Research and Development for Expansion of Radio Resources (JPJ000254)" program.

\ifCLASSOPTIONcaptionsoff
  \newpage
\fi

\bibliographystyle{unsrt}
\bibliography{reference,all_MIC_Literature,export}

\end{document}


%

\appendices
\section{Proof of the First Zonklar Equation}
Appendix one text goes here.

\section{}
Appendix two text goes here.


\ifCLASSOPTIONcaptionsoff
  \newpage
\fi

\begin{IEEEbiography}{Michael Shell}
Biography text here.
\end{IEEEbiography}

\begin{IEEEbiographynophoto}{John Doe}
Biography text here.
\end{IEEEbiographynophoto}


\begin{IEEEbiographynophoto}{Jane Doe}
Biography text here.
\end{IEEEbiographynophoto}




\end{document}